\input amstex
\documentstyle{amsppt}
\NoBlackBoxes
\ifx\undefined\rom
  \define\rom#1{{\rm #1}}
\fi
\ifx\undefined\curraddr
  \def\curraddr#1\endcurraddr{\address {\it Current address\/}: #1\endaddress}
\fi

\topmatter
\title 
The Warwick Automatic Groups Software
\endtitle
\rightheadtext{Warwick Automatic Groups Software}
\author
Derek F. Holt
\endauthor
\address 
Mathematics Institute, University of Warwick,
Coventry CV4 7AL, U.K.
\endaddress
\email dfh\@maths.warwick.ac.uk \endemail

\keywords Automatic Groups, Computer Programs \endkeywords
\subjclass Primary 20F10, 20-04, 68Q40; Secondary 03D40\endsubjclass
\abstract 
This paper provides a description of the algorithms
employed by the Warwick {\smc Automata} package
for calculating the finite state
automata associated with a short-lex automatic group.
The aim is to provide an overview of the whole process, rather than
concentrating on technical details, which have been already been published
elsewhere.  A number of related programs are also described.
\endabstract

\thanks
This paper is in final form and no version of it will be submitted for
publication elsewhere.
\endthanks
\endtopmatter

\document

\head 1. Introduction                
\endhead                             

The purpose of this paper is to provide a description of the algorithms
employed by the Warwick {\smc Automata} package
for calculating the finite state
automata associated with a short-lex automatic group. Such a description,
with full technical details and proofs, has already been published in \cite 2.
Our aim here, however, is to attempt to give a less technical overview of the
full procedure.
This overview also serves to present a more accurate and up-to-date
description of the procedure than \cite 2.
Some recent examples on which {\smc Automata} has been successfully run
will then be discussed.
Finally, we shall describe some associated programs,
which perform various functions on the output. All of these programs are
available from Warwick by anonymous {\bf ftp}. (Please e-mail the author for
details.)

It will be assumed that the reader has some familiarity
with the definition and elementary theory of finite state automata.
The best general reference for the theory of automatic groups is the six-author
book \cite 1.
If $A$ is a finite set (the alphabet), then $A^*$ will denote the set of all
strings of elements of $A$, including the empty string, which will be
denoted by $\varepsilon$. Elements of $A^*$ will be referred to as {\it words}
in $A$.
In this paper, $A$ will always denote a finite set which generates a group
$G$ as a monoid, so each word in $A$ can be regarded as an element of $G$. If
$u$ and $v$ are words in $A$, then we shall use $u = v$ to mean that $u$ and
$v$ are equal as elements of $G$, whereas $u \equiv v$ will mean that $u$ and
$v$ are equal as strings.

Let $A$ and $G$ be as above.
Then $G$ is said to be automatic (with respect to $A$),
if it has an automatic structure.
This consists of  a collection of finite state automata.
The first of these, usually denoted by $W$, is called the {\it word-acceptor}.
It has input alphabet $A$, and accepts at least one word in $A$ mapping
onto each $g \in G$.
The remaining automata $M_x$, are called the {\it multipliers}.
There is one of these for each generator $x \in A$, and also one for $x = Id_G$
(which will be denoted by $M_\$$).
They read pairs of words $(w_1,w_2)$, with $w_1, w_2 \in A^*$,
and accept such a pair
if and only if $w_1$ and $w_2$ are both accepted by $W$,
and $w_1x = w_2$.
(We shall clarify the technical problem of dealing with words $w_1$ and $w_2$
of different lengths at the beginning of the next section.)

Notice that this definition apparently depends on the particular monoid
generating set $A$
for the group $G$, as well as $G$ itself, but it is not too difficult to prove
(see Section 2.4 of \cite 1) that the property of being automatic is preserved
under change of generators, and so we can unambiguously say that a particular
group is automatic or not automatic.
Although it is not necessary theoretically, the programs we are discussing
require that $A$ is closed under inversion  (that is, for each $x \in A$,
there exists $y \in A$ with $xy = Id_G$), and so we shall assume this from
now on.

The definition says nothing specific about the nature of the words accepted
by the word-acceptor. We shall be dealing here exclusively with
so-called {\it short-lex} word-acceptors, which we shall now describe.
For this, we require an
{\it ordered} monoid-generating set $A$ for the group $G$.
For words $u$ and $v$ in $A$, we write $u \prec v$ if either $u$ is shorter than
$v$, or $u$ and $v$ have the same lengths, but $u$ comes earlier than $v$ in
the lexicographical ordering of $A^*$ relative to the given ordering of $A$.
So, for example, if $A = \{a,b\}$ with $a < b$, then
$\varepsilon \prec a \prec b \prec aa \prec ab
\prec ba \prec bb \prec aaa \prec aab$.
This is known as the {\it short-lex} ordering on $A^*$.
A word $w$ in $A$ is said to be {\it irreducible}, if it is the
least word in the short-lex ordering which is equal to $w$ as an element of $G$.
Otherwise, it is called {\it reducible}.

The group $G$ is called {\it short-lex automatic} with respect to the
order monoid-generating set $A$,
if it has a short-lex automatic structure.
This is an automatic structure in which $W$ accepts precisely the
irreducible elements in $A^*$.
In particular, it accepts a unique word for each group element, and accepts
the empty string as representative of $Id_G$.
Note that the equality multiplier
$M_\$$ accepts $(u,v)$ if and only if $u \equiv v$ and $u$ is accepted by $W$.

It turns out that the existence of a short-lex automatic structure 
for $G$ definitely does depend on both the choice and the ordering of $A$.
Some important classes of groups, including word-hyperbolic and finitely
generated abelian groups have short-lex automatic structures for any
choice of ordered generating set.
Other groups, including many Euclidean groups (finitely
generated free abelian-by-finite) have short-lex automatic structures
for some choices but not others.
Finally, there are some groups, such as the braid groups $B_n$ for
$n \geq 4$, which are automatic, but probably have no short-lex automatic
structures at all.

The Warwick {\smc Automata} package can (currently)
only handle short-lex automatic
structures.
It cannot provide any information at all on whether a group might be automatic
if it has no such structure. 
However, if the group does have such a structure using a particular ordered
generating set $A$ then,
in principle, given enough time and space, the programs would be capable of
verifying this fact, and calculating the associated finite state automata
explicitly.
I am not aware of any existing software
which can do the same thing for automatic structures which are not short-lex.
In Chapter 6 of \cite 1, an algorithm is described
for deciding whether a given collection of
finite state automata consists of the automatic structure for some
(not necessarily short-lex) automatic group $G$.
This algorithm has two shortcomings. Firstly, it does not necessarily determine
the group $G$ up to isomorphism (although this would probably be fairly clear
in most cases). Secondly, and more seriously, it does not make any
suggestions as to how to come up with the automata for testing.

The {\smc Automata} package for short-lex automatic groups takes a finite
group presentation as input (it is proved in Theorem 2.3.12 of \cite 1 that
all automatic groups are finitely presented).
The ordering of the generating set is defined by the order in which the
generators are input, and inverses are inserted directly after each generator,
unless they are explicitly specified elsewhere.
The program first attempts to produce
likely candidates for the automata in a short-lex automatic structure for $G$.
It then uses a
method similar to that described in Chapter 6 of \cite 1 to verify the
correctness of these automata. However, in the short-lex situation, the
details are considerably simpler than in the general case and, furthermore,
it is verified that the automata really do correspond to the given group
$G$. Although it is perhaps a somewhat rash claim to make, I believe that
the fact that these programs effectively verify the correctness of the results
that they produce renders the results intrinsically more reliable than those
produced by many computer programs.

The remainder of the paper is organized as follows.
In Section 2, we discuss two-variable finite state automata, logical
operations on finite state automata, and introduce the concept of
word-difference automata, which is central to the main algorithm.
In Section 3, we describe the algorithm itself.
In Section 4, we give a few examples
with time and space requirements, in order to give an idea of the current
scope of the procedure. Finally, in Section 5, we discuss briefly some
related programs that are available with the Warwick package,
which take a short-lex automatic structure as input.
These are really just operations on finite state automata, and are
not inherently dependent on the group structure.
They include methods of determining the cardinality of the
accepted language, enumerating the accepted
language, and determining the growth series of the accepted language.
By applying the last of these to the word-acceptor in a short-lex
automatic structure, we obtain the growth series of the group.

Nearly all of the code in the Warwick package was written by the author,
David Epstein and Sarah Rees. A procedure for minimizing the number of states
in a finite state automaton accepting a given language was written by
Dave Pearson. The growth-series code was written by Uri Zwick.

\head 2. Finite State Automata       
\endhead                             

Let $X$ be a finite state automaton with input alphabet $A$. We denote
the set of states of $X$ by $\Cal{S}(X)$ and the set of accept states by
$\Cal{A}(X)$. The programs in {\smc Automata} deal exclusively with
{\it partial deterministic automata}, and all automata mentioned
in this paper are of this type. (See Definition 1.2.5. of \cite 1.)
They have a unique initial state $\sigma _0(X)$,
and for each $x \in A$ and $\sigma \in \Cal {S}(X)$, there is either a unique
transition from $\sigma$ with label $x$, or no transition. In the first case,
we denote the target of the transition by $\sigma ^x$, and in the second case
we say that $\sigma ^x$ is undefined. We can then define $\sigma ^u$, for
$\sigma \in \Cal {S}(X)$ and $u \in A^*$ in the obvious way. 

The automata in an automatic structure of a group $G$ may be one-variable or
two-variable. The word-acceptor is one-variable, and has input alphabet
a monoid generating set $A$ of $G$. The multipliers (and other automata, to
be introduced later) are two-variable, and the idea is that they should read
two words $u$ and $v$ in $A$ simultaneously, and at the same rate. This
creates a technical problem if $u$ and $v$ do not have the same length. To
get round this, we introduce an extra symbol  $\$$, which maps onto the
identity element of $G$, and let $A^\dagger = A \cup \{\$\}$. Then if
$(u,v)$ is an ordered pair of words in $A$, we adjoin a sequence of
$\$$'s to the end of the shorter of $u$ and $v$ if necessary, to make them
both have the same length. The resulting pair will be denoted by
$(u,v)^\dagger$, and can be regarded as
an element of $(A^\dagger \times A^\dagger)^*$. Such
a word has the property that the symbol $\$$ occurs in at most one of $u$
and $v$, and only at the end of that word, and it is known as a {\it padded}
word in $A \times A$. We shall assume from now on,
without further comment, that all of the two-variable automata that arise
have input language $A^\dagger \times A^\dagger$ and accept only padded words.

The {\smc Automata} package makes heavy use of logical operations on
finite state automata.
Let $X$ and $Y$ be two finite state automata, with the same input
alphabet, and let the accepted languages be $L(X)$ and $L(Y)$, respectively.
Then there are straightforward algorithms for constructing finite state
automata $X\cdot Y$,
$X^\prime$, $X \wedge Y$ and $X \vee Y$, with accepted languages
$L(X)L(Y)$, $A^* \setminus L(X)$, $L(X) \cap L(Y)$
and $L(X) \cup L(Y)$, respectively.
If $Z$ is a two-variable automaton accepting padded words in  $A \times A$
then there is one-variable automata $E(Z)$, with input alphabet
$A$, that accepts $u \in A^*$ if and only if there exists $v \in A^*$ with
$(u,v)^\dagger \in L(Z)$. (Clearly the same applies with $u$ and $v$
interchanged, and we can construct the corresponding ``for all'' automata by
using a combination of negation and ``there exists''.) These operations are
all implemented in {\smc Automata}. The algorithms are described in Chapter
1 of \cite 1, and we shall not go into details here, except to mention that
the construction of $E(Z)$ amounts to constructing a deterministic
automaton equivalent to a non-deterministic one, and potentially involves
increasing the number $n$ of states of $Z$ to $2^n$ states of $E(Z)$.
This is one of the
reasons why the programs can be very space-demanding for larger examples.
{\smc Automata} also contains a procedure for replacing an automaton $X$
with an equivalent one (that is, one accepting the same language) with the
minimal number of states. This is applied routinely after all of the
logical operations that we have just described.

The principal difference in emphasis between the current paper and \cite 2
is that we shall describe as much of the algorithm as possible in terms
of these logical operations on finite state automata. In \cite 2, we
gave more detailed direct descriptions of the individual procedures,
and put much more emphasis on the string-rewriting aspects.
I hope that the current approach is conceptually easier to understand.

We shall now introduce the concept of {\it word-difference automata}, which
play a central role in the programs.
Assume, as usual, that $A$ is a monoid generating set for a group $G$. Then
a two-variable automaton $Z$ accepting padded words in $A \times A$ is
called a word-difference automaton, if there is a function
$\alpha:\Cal{S}(Z) \rightarrow G$ such that
\roster
\item"(i)" $\alpha(\sigma_0(Z)) = Id_G$, and
\item"(ii)"  for all $x,y \in A^\dagger$
and $\tau \in \Cal{S}(Z)$ such that $\tau ^{(x,y)}$ is defined, we have
$\alpha(\tau ^{(x,y)}) = x^{-1}\alpha(\tau)y$.
\endroster

We shall assume that all states $\tau$ in a word-difference automaton $Z$ are
accessible; that is, there exist words $u,v$ in $A$ such that
$\sigma _0(Z)^{(u,v)^\dagger} = \tau$.  It follows from
Properties (i) and (ii) that $\alpha(\tau) = u^{-1}v$, and so the map
$\alpha$ is determined by the transition function of $Z$.
We shall also assume that each accept state is mapped by $\alpha$ onto
the same element $g$ of $G$. We then say that $Z$ is a word-difference
automaton accepting $g$.

It is proved in Theorem 2.3.2 of \cite 1 and Theorem 2.3 of \cite 2 that
any automatic structure has the so-called {\it fellow-traveler} property.
This means that there is a finite set $\Cal {WD}$ (the word-differences) of
elements of $G$ with the following property.
If $u,v \in L(W)$ and either $u = v$ or $u^{-1}v = x$
for some $x \in A$, and $p(u)$ and $p(v)$ are prefixes of $u$ and $v$
having the same length, then $p(u)^{-1}p(v) \in \Cal {WD}$.
We also insist that $x \in \Cal {WD}$ for all  $x \in A \cup \{Id_G\} $.
We can therefore define associated word-difference automata $WD_x$,
for $x \in A \cup \{Id_G\} $, in which the set of states is $\Cal {WD}$,
$\alpha$ is the identity map, and $\tau ^{(x,y)}$ is defined if and only if
$x^{-1}\tau y \in \Cal {WD}$, in which case it is equal to $x^{-1}\tau y$.
The initial state is $Id_G$, and there is a single accept state $x$.

It is clear from all of this (c.f. Theorem 2.3.4 of \cite 1
and Corollary 2.4 of \cite 2)
that the multipliers $M_x$ in an automatic structure
accept $(u,v)^\dagger$ if and only if $WD_x$ accepts $(u,v)^\dagger$
and the word-acceptor $W$ accepts both $u$ and $v$.
As we shall see in the next section, this fact can be used to construct them.

There are two other word-difference automata which arise in the algorithms.
We define $L_1$ to be the set of padded words $(u,v)^\dagger$ in $A \times A$
with the following properties
\roster
\item"(i)" $u = v$ in $G$;
\item"(ii)" $v$ is irreducible;
\item"(iii)" $u$ is reducible;
\item"(iv)" $u = tx$ with $t \in A^*, x \in A$, and $t$ is irreducible;
\item"(v)" $u = ys$ with $s \in A^*, y \in A$, and $s$ is irreducible.
\endroster
So, if we think of the elements of $L_1$ as reduction rules for $G$, where
we replace occurrences of the left hand sides of rules by the corresponding
right hand sides, then $L_1$ is the minimal set of reduction rules necessary
to reduce any word in $A$ to its irreducible representative.
If $G$ has a short-lex automatic structure with respect to $A$, then
it follows from the fellow-traveler property that the automaton
$WD_\$$ accepts all words  $(u,v)^\dagger$ in $L_1$.
(This is because $WD_x$ accepts $(t,v)^\dagger$ with $t$ as in (iv).)
Since the word-acceptor $W$ accepts
precisely the set of irreducible words, it is clear that the words
$(u,v)^\dagger$ satisfying properties (ii)--(v) can be recognized by some
finite state automaton $X$.
Since the words accepted by $WD_\$$ satisfy property (i),
it follows that $X \wedge WD_\$ $ is a word-difference
automaton accepting the language $L_1$.
More generally, we shall say that a word-difference automaton
accepting $Id_G$ has property $\Cal {P}_1$, if its language contains $L_1$.

We define $L_2$ to be the same as $L_1$, except that we omit property (v)
from the definition; so $L_2$ is usually larger than $L_1$.
Again, there is a word-difference automaton accepting $L_2$, and
a word-difference automaton accepting $Id_G$
is said to have property $\Cal {P}_2$ if its language contains $L_2$.
In particular, $WD_\$$ has property $\Cal {P}_2$.

In the programs, we shall be constructing automata $D_1$ and $D_2$ which
we hope will have properties $\Cal {P}_1$ and $\Cal {P}_2$, respectively.
(These are referred to as the {\it principal transition table} and the
{\it product transition table} in \cite 2.)

\head 3. Description of the Algorithm                
\endhead                             

To understand the algorithm, it is essential to distinguish between
the theoretical and the practical points of view. The material in
Section 2 was completely theoretical. We assumed the existence of
a short-lex automatic structure for a group $G$, and used the
associated automata $W$ and $M_x$ to deduce the existence of $WD_x$,
and word-difference automata with properties $\Cal {P}_1$ and $\Cal {P}_2$.

The algorithms also assume the existence of an automatic structure,
but their aim is to construct this structure,
and $W$ and $M_x$ are not explicitly available at the beginning.
Thus we cannot use them to construct the word-difference
automata $D_1$ and $D_2$.
On the contrary, we first guess candidates for $D_1$ and $D_2$,
and then use these to construct candidates for $W$ and $M_x$,
and finally check the correctness of the complete automatic structure.

Another point to be borne in mind is that we cannot solve the word problem
in $G$ until the whole process is complete. This means that we cannot
necessarily calculate the map $\alpha$ in the word-difference automata
precisely, even though we know theoretically that it exists. The best
we can do is to calculate a map $\hat \alpha$ from the set of states to $A^*$,
where $\hat \alpha(\sigma)$ is the least element of $A^*$ (under the short-lex
ordering) which is currently known to be equal to  $\alpha(\sigma)$ in $G$.
In the construction of $D_1$ and $D_2$, we really want the states to be
elements of $\Cal {WD}$ and $\alpha$ to be the identity map, and so in
practice, the states are elements of $A^*$ which we cannot currently
reduce, and $\hat \alpha(\sigma)$ is the identity map. Whenever we
discover that two states are equal in $G$, we amalgamate them.

We can now give a somewhat simplified breakdown of the complete
algorithm. The input is an ordered monoid-generating set $A$ of a finitely
presented group $G$, together with defining relations for $G$ given as
words in $A$. In principal, it will succeed if and only if $G$ is
short-lex automatic with respect to $A$, in which case it will construct
the finite state automata of the associated automatic structure.
(After minimization, these automata are essentially uniquely determined,
because their languages are uniquely determined by $A$ and $G$.)

\roster
\item"{\bf Step 1}." Construct a word-difference automaton $D_1$
accepting the identity, which we hope has property $\Cal {P}_1$.
\item"{\bf Step 2}." Construct a candidate for the word-acceptor $W$.
Construct the equality multiplier $M_\$$, which accepts  the
padded word $(u,v)^\dagger$ if and only if $u=v \in L(W)$.
\item"{\bf Step 3}." Construct a word-difference automaton $D_2$
accepting the identity, which we hope has property $\Cal {P}_2$.
\item"{\bf Step 4}." Construct candidates for the multipliers $M_x$, for all
$x \in A$.
If, during this construction, we find distinct  words
$u$ and $v$ in $A$, both accepted by $W$, with $u = v$ in $G$, then return
to Step 1.
\item"{\bf Step 5}." Check correctness of $D_2$. If incorrect, return to Step 3.
\item"{\bf Step 6}." Check correctness of $W$ and $M_x$. If incorrect, return to
Step 1.
\endroster

We shall describe these steps in more detail.
All of them except for Step 1 are essentially logical operations on
finite state automata, and so we leave Step 1 until last.
Assume then that we have completed Step 1, and
let $GT$ be the two-variable automaton which, for $u,v \in A^*$ accepts
$(u,v)^\dagger$ if and only if $u \succ  v$.
In Step 2, we construct $W$
using the standard finite state automaton operation
$$ W = (A^*\cdot E(D_1 \wedge GT)\cdot A^*)^\prime. $$
Expressed in words, $W$ accepts $w \in A^*$ if and only if it does not have
a substring $u$ such that there exists $v \in A^*$ with
$u \succ v$ and $(u,v)^\dagger$ accepted by $D_1$.
It is clear that, if $D_1$ really does have property $\Cal {P}_1$, then
$W$ will accept the set of irreducible words; in other words, it will be
correct.

The converse is not necessarily true; that is, it is possible for $W$ to
be correct and for $D_1$ not to accept everything in $L_1$. However,
what is true is that if $W$ is correct and $u$ is reducible (so not
accepted by $W$), then there is a substring $s$ of $u$ such that $D_1$
accepts $(s,t)^\dagger$ for some $s \succ t$, and it is not difficult
to use $D_1$ to find $s$ and $t$. (This technique is contained in the
proof of Theorem 4.1 of \cite 2, and we shall not give further details
here.) Thus, if $W$ is correct, then we can effectively use $D_1$ to
reduce words in $A^*$ to their irreducible representatives. That is, we
can solve the word problem in $G$ and calculate the maps $\alpha$
associated with the word-difference automata correctly.
The subsequent steps of the algorithm use this technique under the
assumption that $W$ is indeed correct.
As we shall see shortly,
if this is the case, and $G$ really is short-lex automatic with respect to $A$,
then it is guaranteed that after a finite number of iterations through the
loop in Steps 3--5, the multipliers $M_x$ will also be correct, and so Step 6
will succeed. 

Another significant property of the $W$ constructed is that, even if it is
incorrect, it always accepts at least one word for each group element.
To see this, let $g \in G$ and let $u$ be a word in $A$ mapping onto $G$.
Then we can use $D_1$ to reduce $u$ to a word $v$ accepted by $W$, and since
we are assuming that $D_1$ accepts the identity, we have $u = v = g$.
This means that the language accepted by $W$ always contains the irreducible
words; that is, it can be too large but never too small.
Since Step 6 effectively checks that $W$ accepts only irreducible words,
we know that, if Step 6 succeeds, then the automatic structure that we have
constructed is for the given group $G$. In the more general situation
described in Chapter 6 of \cite 1, nothing about $W$ is assumed, and so it
is possible, even if the verification process succeeds, that the automatic
structure is for a proper factor group of the given group $G$.

The equality multiplier $M_\$$ constructed in Step 2 is not really
necessary but it is useful to have it available
in Step 6.

In Step 3, our aim is to make $D_2$ equal to the automaton $WD_\$$
defined in Section 2. For this, we need a candidate for the word-difference
set $\Cal {WD}$. For our initial candidate (on the first call of
Step 3), we let  $\Cal {I}$ be the set of all images $\hat \alpha(\tau)$ for
$\tau \in \Cal {S}(D_1)$, and put
$$\Cal {S}(D_2) = \Cal {I} \cup \Cal {I}^{-1} \cup A \cup \{\varepsilon\}.$$
We do this because we know that $\Cal {WD}$ is closed under inversion
(since $A$ is), and we require it to contain $A\cup \{\varepsilon\}$.
We then define $D_2$ exactly as specified for $WD_\$$ in Section 2.
That is, for $\tau \in \Cal {S}(D_2)$ and $x,y \in A^\dagger$,
$\tau ^{(x,y)}$ is defined if and only if
$x^{-1}\tau y \in \Cal {S}(D_2)$, in which case it is equal to $x^{-1}\tau y$.
The initial state and the only accept state is $\varepsilon$.
The expressions $x^{-1}\tau y$ are calculated using the word reduction
algorithm mentioned above.
Remember that, if $W$ is correct, then so is the word reduction algorithm,
and so $\Cal {S}(D_2)$ can be regarded as a subset of $G$ and $\hat \alpha =
\alpha$.
It follows that, if $W$ is correct and
$\Cal {S}(D_2)$ really does contain $\Cal {WD}$, then
$D_2$, as constructed, will have property $\Cal {P}_2$.
However, in many examples $\Cal {S}(D_2)$ will not contain $\Cal {WD}$
on the first iteration, and Step 5 tests for this, and produces some explicit
elements in $\Cal {WD} \setminus \Cal {S}(D_2)$ if it is nonempty.
If $G$ is
short-lex automatic with respect to $A$, then $\Cal {WD}$ is finite, and so
we will only need to go through this iteration a finite number of times
(assuming of course that $W$ is correct).

We next describe Step 4.
We construct the multipliers exactly according to the definition of
{\it standard automata} in Definition 2.3.3 of \cite 1. For each $M_x$, the
set of states is $\Cal {S}(W) \times \Cal {S}(W) \times \Cal {S}(D_2)$
with initial state $(\sigma _0(W),\sigma _0(W),\sigma _0(D_2))$,
and for $x,y \in A^\dagger$ and $(\sigma _1,\sigma _2, \tau) \in \Cal {S}(M_x)$,
the transition $(\sigma _1,\sigma _2, \tau)^{(x,y)}$ is defined to be equal to
$(\sigma _1^x,\sigma _2^y,\tau^{(x,y)})$ provided all three components are
defined, and undefined otherwise. The $M_x$  for different $x$
differ from each other only in their accept states; the accept states
of $M_x$ are all $(\sigma _1,\sigma _2, \tau)$ such that
$\sigma _1, \sigma _2 \in \Cal {A}(W)$ and $\hat \alpha (\tau) \equiv x$. 
However, once the $M_x$ have been constructed in this fashion, they are
immediately minimized, after which they no longer have the same sets of
states, and are no longer necessarily word-difference automata.

It sometimes happens that, during the course of the construction in Step 4,
we find distinct words $u,v$ in $A$, both accepted by $W$, for
which $u = v$ in $G$. This means that $W$ is accepting two words for the
same group element, and so it cannot be correct. We can therefore return
immediately to Step 1.
(This happens when there is $(u,v) \in L(D_2) \setminus L(D_1)$, where
$(u,v)$ ought to be in $L(D_1)$, and so $D_1$ is incorrect.)
Although we know of no theoretical reason why this
should be the case, we have never yet encountered an example in which
we have got through to Step 6 and the axiom checking has failed. It seems
that, if our automata are wrong, then this fact becomes clear either at
this stage or at Step 5.

Note that, if $W$ is correct (which ensures that $\Cal {S}(D_2)$ is a subset
of $G$ and $\alpha$ is the identity map) and $D_2$ accepts everything in $L_2$
(defined in Section 2),
then $\Cal {S}(D_2)$ must contain $\Cal {WD}$, and so the multiplers $M_x$
will be correct. Suppose  that $W$ is correct but $M_x$ is not correct, for
some $x \in A$. From the construction of $M_x$, it is clear that, if
$M_x$ accepts $(u,v)^\dagger$, then $W$ accepts $u$ and $v$, and
$ux = v$ in $G$. Thus, there exist words $u$ and $v$ in $A^*$ such that
$W$ accepts $u$ and $v$ and $ux = v$ in $G$, but $M_x$ does not accept
$(u,v)^\dagger$. But then, since $v$ is the unique word accepted by $W$
that is equal to $ux$ in $G$, $M_x$ cannot accept $(u,w)^\dagger$ for
any $w \in A^*$. So we can test for this, by constructing the automata
$W \wedge E(M_x)'$ for each $x \in A$, and we do this in Step 5.
If all of these automata have empty accepted language and $W$ is correct,
then the $M_x$ must also be correct, and we proceed to Step 6 for the final
verification. Otherwise, we find some explicit words $u$ which are
not accepted by some $W \wedge E(M_x)'$, and we let $v$ be the reduction
of $ux$ (using $D_1$); then $W$ accepts $v$ and $ux = v$ in $G$.
Since $M_x$ does not accept $(u,v)$, there must be prefixes $p(u)$ and
$p(v)$ of $u$ and $v$ of the same length, such that $p(u)^{-1}p(v) \notin
\Cal {S}(D_2)$. We therefore return to Step 3, and adjoin all such
elements and their inverses to $\Cal {S}(D_2)$, and then recalculate
$D_2$ itself.

Step 6 depends on the following logical automaton operation. If $Z_1$
and $Z_2$ are any two two-variable automata with the same input alphabet
$A^\dagger \times A^\dagger$, then we define their composite to be the
two-variable automaton with language
$$ \{ (u,v)^\dagger  | \exists w \in A^* : (u,w)^\dagger \in L(Z_1)
                                \text{ and } (w,v)^\dagger \in L(Z_2) \}.$$
In particular, for multipliers $M_x$ and $M_y$, we denote their composite by
$M_{xy}$. The composite operation is easily seen to be associative, and
so we can define $M_u$ for any $u \in A^*$ by repeated application, where
$M_\varepsilon = M_\$ $. It is also easy to see that, if $M_u$ accepts
$(v,w)^\dagger$, then $vu=w$ in $G$. The axiom checking process essentially
verifies the converse of this statement. More precisely, let $\Cal {R}$ be
the set of defining relators of $G$, as a monoid (so it includes the
inverse relators $xx^{-1}$, for all $x \in A$). Then, it is proved in
Theorem 2.5 of \cite 2, that if $M_r$ has the same language as $M_\$$ for
all $r \in \Cal {R}$, then the automata $W$ and $M_x$ are correct.
So we can check this by repeated application of the composite operation.
Since this operation involves making a non-deterministic automaton
deterministic, it potentially involves an exponential increase in the
number of states. We try to keep this under control by minimizing after
every operation. With a long relator $r$ it seems to be most efficient
to split $r$ roughly in half as $r \equiv st$ and then form $M_s$ and
$M_t$, and construct $M_r$ as the composite of these. 
In most examples, Step 6 takes the longest time, and is often the most
space-demanding.

Finally we turn to Step 1. Our current approach to this is to use the
Knuth-Bendix procedure on the group presentation to derive a list of
reduction rules of the form $(u,v)$, where $u$ and $v$ are words in $A$,
$u = v$ in $G$ and $u \succ v$. Since we devoted a large part of \cite 2
to the Knuth-Bendix procedure, we will not discuss it in detail here, but
concentrate on its application to Step 1 of the procedure.
At any time, the existing rules are used to reduce words $w$ in $A$ to their
least possible representative, by substituting $v$ for $u$ whenever $u$ is
a substring of $w$ and $(u,v)$ is a rule. Occasionally, the Knuth-Bendix
procedure terminates with a finite complete set of reduction rules,
in which case these rules suffice to reduce any word to its irreducible
representative in $G$. But this does not usually happen.
For any particular $w$, it is guaranteed that we will eventually generate
enough rules to correctly reduce $w$, but we have no way in general of
knowing how long this will take, or even whether we can currently reduce
$w$ correctly.
(For example, the group itself could be trivial, but it could take a long
time for this fact to emerge.)

At regular intervals, we interrupt the Knuth-Bendix procedure, and form a 
candidate for the word-difference automaton $D_1$ using the current set
of rules. We do this as follows.  $\Cal {S}(D_1)$ is initially empty.
Consider each rule $(u,v)$ in turn, and
let $u = x_1x_2 \ldots x_m$ and $v = y_1y_2 \ldots y_n$ with $m \ge n$.
We pad $v$ if necessary in the usual way, by putting $y_i = \$ $ for
$n < i \le m$. Let $d_0 = \varepsilon$  and for $1 \le i \le m$ let
$d_i$ be the reduction, using the current rules, of
$x_i^{-1}d_{i-1}y_i$ (so $d_i$ is a word-difference).
Then we adjoin the states $d_i$ to $\Cal {S}(D_1)$ and the transitions
$d_{i-1}^{(x_i,y_i)} = d_i$ to $D_1$ if any of them are not there already.
Of course $\alpha$ maps a state to its image in $G$, and the start and accept
states are both $\varepsilon$.

Suppose $G$ really is short-lex automatic with respect to $A$. Then, as we
saw in Section 2, there is a word-difference automaton that accepts
the language $L_1$.
Since this automaton contains only finitely many states and transitions, it
follows that these states and transitions will be manifested by a finite
subset of $L_1$. Thus once we have run the Knuth-Bendix procedure long
enough to generate all of the rules in this finite subset, and we have
enough rules to correctly reduce each of the transitions $x_i^{-1}d_{i-1}y_i$
described above, the candidate for $D_1$ that we construct will have
the required property $\Cal {P}_1$, and then the remaining steps 2--6 of
the complete procedure will also complete successfully.

The big problem is to recognize when this has happened and to proceed to
Step 2. Currently, we have no completely satisfactory heuristics for this.
Our default method is to wait until several successive candidates for
$D_1$ have been identical; that is to wait until $D_1$ has apparently
become constant. This approach is successful for many straightforward
examples. However, because the current system of rules at any time is
not adequate to reduce all strings correctly, it happens sometimes that
the $D_1$ constructed has superfluous states and accepts strings that are
not in $L_1$. This, in itself, does not matter provided that there are not too
many such superfluous states, but it has the effect that it makes it difficult
to recognize when $D_1$ has property $\Cal{P}_1$, because it never appears to
become stable. In such cases, the user has to decide interactively when to stop.
Of course, if $D_1$ appears to be increasing steadily and rapidly in size,
then this is good evidence that $G$ is not short-lex automatic with respect
to $A$.

Currently, if we need to return to Step 1 after Steps 3 or 6, then we begin
again completely, because it seems to difficult (and space consuming) to attempt
to restart from where we left off last time. As we mentioned above, we have
never encountered an example in which Step 6 failed.
It would probably be a good idea to insert the discovered identity $u = v$
as a new rule if we return to Step 1 after Step 3, but we are not currently
doing that.
Alternatively, we could simply recalculate $D_1$ after adjoining the new
rule $(u,v)$ and not use Knuth-Bendix again at all.
We intend to experiment with such alternative strategies in the near future.
Another simple policy, which we do use sometimes, is to use $D_2$ in place
of $D_1$, which seems reasonable theoretically, since the language of
$D_2$ is usually larger than that of $D_1$. In fact, this is often effective;
the only disadvantage is that it usually increases the space requirements of
Step 2 quite considerably, and so it is not always practical.

As we have seen, Steps 2--6 are more routine and predictable than Step 1,
but they are subject to exponential increase in storage demands at certain
points. Our experience with examples indicates that when $D_1$ is correct
(i.e. it has property $\Cal {P}_1$, then the space requirements will be large,
but not impracticably so, but when $D_1$ is incorrect, the exponential
increases will really occur and cripple the whole process. There is no
known theoretical reason why this should be the case, but it is seems reasonable
to suppose that when $D_1$ is correct then the subsequent calculations are
constructing natural objects associated with $G$, whereas when $D_1$ is
incorrect, the calculations are in some sense meaningless and therefore
random.

This situation indicates that the most important problem now is to try to
find improved methods of carrying out Step 1.  Improvements to the remaining
steps would also be useful of course, but they are more likely to be purely
computational in nature rather than algorithmic. In the Knuth-Bendix procedure,
there is considerable freedom in the choice order in which we process overlaps.
A strategy needs to be found which produces the finite subset of $L_1$ mentioned
above as quickly as possible. There may also be completely different approaches
to Step 1, which do not involve Knuth-Bendix, but none has been found to date.

\head 4.  Examples                   
\endhead                             

A number of relatively easy and straightforward examples in which
the {\smc Automata} program runs successfully were described in \cite 2.
Here we describe some examples of varying difficulty which have arisen
in some natural context during the past few years, for which we have
successfully computed the automatic structure. There are two shell-scripts
available which run the complete procedure described in Section 3, by
running the individual programs that perform the steps in the procedure
in the correct order. These are called {\it automata} and {\it automatac}.
The second of these has often proved more efficacious on the harder examples.

The first example is a group that was used by Charlie Gunn at the Geometry
Center, Minneapolis, in connection with drawing pictures of
tessellations of hyperbolic space with a dodecahedral fundamental region.
Such pictures were used in the film "Not Knot" produced by the Center.
The word-acceptor of the group is used to generate each group-element
exactly once, which results in faster and more accurate drawings.
The presentation is
$$ G_1 = \langle  a,b,c,d,e,f | 
a^4, b^4, c^4, d^4, e^4, f^4,
abAe, bcBf, cdCa, deDb, efEc, faFd \rangle.
$$
(In all of the examples in this section, we use the inverse by case-change
convention; that is $A$ is the inverse of $a$, etc.)
The {\it automatac} program completed in a total of 545 seconds
on a Silicon Graphics Iris XS24 machine, with 32 Megabytes of core memory
available.  The word-acceptor has 48 states and the multipliers range from
145--175 states. This example went twice round the first loop (Steps 1--4) of
the algorithm, and six times round the second loop (Steps 3--6).

In the above example, and in many of those that follow, it would be
possible to eliminate several of the generators before starting.
Experience has shown that this is not a good idea in general.
The programs appear to run best on the natural presentations of the groups,
rather than on presentations that have been manipulated by the user.

The second example is a group arising from a $(2,2,2,2,2,2)$-orbifold
proposed by Bill Thurston. The presentation is
$$ \align G_2 = \langle & a,b,c,d,e,f,g,h,i,j,k,l |  \\
&dcbadc, (ahIE)^2, (beJF)^2,
(cfKG)^2, (dgLH)^2, (ilKJ)^2 \rangle.
\endalign
$$
Initially this caused problems because of the large number of generators,
which meant that the alphabet for the two-variable automata had size
361. To cope with this, we had to provide an optional new method of storage
for automata, essentially storing the transition table as a sparse matrix,
rather than as a simple two-dimensional array. The example then became
tractable. The word-acceptor has 131 states and the multipliers range from
132--232 states. The total time on the SGI machine was 930 seconds.

The Fibonacci groups have been favorite test examples in combinatorial group
theory for many years. The groups $F(2,n)$ (for $n >2$, say) have the
presentation
$$  F(2,n) = \langle a_1, \ldots, a_n | 
a_1a_2=a_3, a_2a_3=a_4, \ldots , a_{n-1}a_n=a_1, a_na_1=a_2 \rangle.
$$
The group is finite when $n=3,4,5$ or $7$. For $n=6$ it is Euclidean
(free abelian-by finite), and is an easy example for {\it automata}.
For $n \geq 8$ and $n$ even, they are known to be hyperbolic, and hence
short-lex automatic. For $n=8$, the $automata$ program ran out of space on
the SGI machine while computing the multipliers. It completed (in about
9 hours cpu-time) on a much slower Solbourne machine with 128 Megabytes
of memory. The word-acceptor has 212 states and the multipliers all have 1861
states. In this example, the {\it automatac} program went twice round the
first loop, then twice round the second loop, and then needed to go another
time around the first loop before completing. In other words $D_1$ was
wrong the first time (with 51 word-differences), then $D_2$ was wrong,
and then it transpired that $D_1$ was wrong again (with 54 word-differences).
There are in fact 58 states (word-differences) in the correct $D_1$
automaton.

We have also successfully computed the automatic structure of $F(2,10)$,
which is larger (and more difficult) than that of $F(2,8)$, but not
disproportionately so. The examples for $n$ odd are much more difficult,
and it is not known whether or not they are hyperbolic.
Since writing the first draft of this paper, the author has rewritten some
of the procedures in {\smc Automata}, using improved data structures (but
no essentially new algorithms), and it has become possible to complete some
larger examples. In particular, $F(2,9)$ has been proved to be short-lex
automatic using its natural presentation, and this is the largest automatic
structure that we have successfully computed to date. The word-acceptor has
3251 states, and the multipliers about 25000 states,
but they had about 860000 states before minimization.
The complete computation took about 12 hours of cpu-time
and used up to about 130 Megabytes of memory on a new SPARCstation 20.
More details of this calculation can be found in \cite 3.
This group had already been proved to be infinite by M. F. Newman in \cite 4;
the automatic structure easily provides an independent proof of this fact,
together with additional information, such as the fact that the generators
have infinite order, which had not been known previously.

We conclude by mentioning some further large examples which we have successfully
computed recently.
The first of these is a presentation for the Picard group $P=SL(2,R)$,
where $R$ is the ring of Gaussian integers, taken from \cite 6:
$$ \align P = \langle & a,t,u,l,j | \\
& a^2 = l^2 = (at)^3 = (al)^2 = (tl)^2 = (ul)^2 = (ual)^3 = j,
j^2=[t,u]=1,  j \text{ central}  \rangle.
\endalign
$$
The word-acceptor has
403 states and the multipliers range from 403--3718.  (The smallest is for the
central generator $j$, which is of course $-I$ in $SL(2,R)$.) The most
remarkable thing about this example is that it was necessary to go 12
times around the second loop in the procedure. This suggests that techniques
should be sought to find more missing word-differences at the same time.

The following presentation was proposed by Andr\'e Rocha:
$$  G = \langle  a,b,c,d,e,f | 
 BceFa, CdfBa, DebCa, EfcDa, FbdEa, ecfdb \rangle.
$$
(This is the fundamental group of a negatively closed Seifert-Weber
manifold.)
The word-acceptor has 1429 states, and the multipliers about 14000 states.

The final presentation was proposed several years ago by Thurston, and 
represents a knot with surgery. (The first two relations define the group
of the knot, and the final relation the surgery. The computation for
the knot group itself, obtained by omitting the final relation, is relatively
straightforward.)
$$  G = \langle  x,y,t | 
 txT = xyx, tyT = xy, (XYxy)^2 = t^3 \rangle.
$$
The word-acceptor has 1188 states, and the multipliers about 9300 states.

\head 5.  Related Algorithms                   
\endhead                             

In this final section, we mention some related algorithms that have been
implemented at Warwick, which take an automatic structure of a group as
input. The first few of these are fairly easy, and the implementations
are distributed with {\smc Automata}. There is a program to enumerate the
accepted words of a finite state automaton up to a given length,
using a depth-first search. When applied to the word-acceptor of a short-lex
automatic group, this produces a unique word for each group element, and is
the basis of many of the applications of the programs to hyperbolic geometry,
particularly to drawing pictures. There is a fast and simple program to
decide whether the language of an automaton is finite and, if so, to
determine its order. This can be used on the word-acceptor to determine
the order of a finite group. (For most presentations of finite groups,
Todd-Coxeter coset enumeration provides the fastest mechanical way of
determining the order, but there are some cases in which the automatic group
method works faster.)

If a group is word-hyperbolic, then the geodesics between any two vertices
in the Cayley graph have the fellow-traveler property - that is they remain
within a bounded distance of each other. It has recently been proved by
Panos Papasoglu (see \cite 5) that the converse of this statement is also true;
if the geodesics in the Cayley graph have the fellow-traveler property,
then the group is hyperbolic.
This means that there should be an automatic structure
for the group in which the word-acceptor accepts all shortest words for each
group element. There is a program available  with the package
which attempts to calculate this word-acceptor. If it terminates, then the
geodesic word-acceptor that is calculated is guaranteed to be correct and so,
by the result of Papasoglu, we have verified that the group is hyperbolic.

The algorithm used is as follows. 
We assume that some short-lex automatic structure for the group $G$ has already
been calculated.
Assume also that we have a list $\Cal {WDG} $ of some or all of the
word-differences that occur between geodesics in the Cayley graph of $G$,
and that $\Cal {WDG} $ contains the set $\Cal {WD} $ (see Section 3) that
was calculated with the short-lex automatic structure.
Let $WDG$ be the corresponding two-variable word-difference automaton
which accepts $(u,v)$ if and only if $u=v$, $u$ and $v$ have equal lengths,
and the word-differences arising from the two paths $u$ and $v$
all lie in $\Cal {WDG}$.
Now calculate a sequence $W_i (i \geq 0)$ of finite state automata,
where $W_0$ is the short-lex word-acceptor, which we already have,
and 
$$ W_{i+1} = \{ u | \exists v \in L(W_i) \text{ with } (u,v) \in L(WDG) \}.$$
Then $L(W_i) \subseteq L(W_{i+1})$ for all $i$, and 
if $L(W_i) = L(W_{i+1})$ for some $i$, then  $L(W_i)$ is the set of all
geodesic words, and so $W_i$ is the required geodesic word-acceptor.
(To see this, let $u$ be any geodesic word, let $v$ be the short-lex reduced
word equal to $u$ in $G$, and let $u=u_0,u_1, \ldots, u_n=v$ be a reduction
of $u$ to $v$ where, for each $i$, $(u_i,u_{i+1})$ is in the language $L_2$
defined in Section 2. Since $u$ is a geodesic, the $u_i$ all have the same
length and, since we are assuming that $\Cal {WDG} $ contains $\Cal {WD} $,
we have $(u_i,u_{i+1}) \in L(WDG)$ for all $i$, and so $u \in L(W_n)$.)
So if the algorithm terminates, then we know that the answer is correct.
Conversely, if $\Cal {WDG}$ contains all the word-differences that arise
from pairs $(u,v)$, where $u$ is a geodesic and $v$ is the short-lex reduced
word equal to $u$ in $G$, then the algorithm will terminate; in fact $W_1$
will already be the correct geodesic word-acceptor.
The only problem with all of this is how to find the set $\Cal {WDG} $.
We can try $\Cal {WDG} = \Cal{WD}$ as a first attempt,
but that is often not large enough. 
Currently, our method is to choose random geodesics $u$ up to a specified length
(40 or 50 is usually adequate in the examples that we are considering),
reduce them to their short-lex minimal equivalents $v$, calculate the
word-differences that arise from $(u,v)$, and adjoin them to $\Cal {WDG}$
if they are not there already. A deterministic method for doing this has been
proposed by David Epstein but has not yet been implemented.

For example, in the group $G_1$ of Section 4, $| \Cal{WD} | = 75$ and
the algorithm does not appear to terminate with $\Cal {WDG} = \Cal{WD}$, but
by choosing random geodesics we rapidly increase $\Cal{WDG}$ to size 103,
and the procedure then terminates with the geodesic word-acceptor $W_1$ 
having 64 states.
In $G_2$, $\Cal{WD}$ has size 99 and is already adequate to calculate
the geodesic word-acceptor, which has 156 states.
We have not yet succeeded in calculating the geodesic word-acceptor
for $F(2,8)$.

Finally, we should like to mention the growth-rate program, written by
Uri Zwick. This is not supplied with {\smc Automata}, but a binary
executable for Sun 4 architecture is available on request.
The growth-function for a group $G$ with a given monoid generating set $A$ is
defined to be the formal power series
$f_{G,A}(t) = \sum_{i \geq 0} w(i)t^i$, where $w(i)$ is the number of elements
of $G$ for which the shortest word in $A$ has length $i$.
If we have a short-lex word-acceptor $W$ for $G$ with respect to $A$, then
clearly $w(i)$ is equal to the number of words of length $i$ in $L(W)$,
and so $f_{G,A}(t)$ can be calculated from $W$. It can be proved that the
generating function of the language of a finite state automaton is always
a rational function of $t$, and the result is output in this form by Zwick's
program. For example, the growth function for the free group on two
generators is $(1+t)/(1-3t)$, whereas the growth-rate for the example $G_1$
in Section 4 is $( 1 + 3 t + 3 t^2 +  t^3 ) /( 1 - 9 t + 9 t^2 -  t^3 )$.

\enddocument